\documentclass[12pt]{l4dc2022}

\usepackage{amsmath,amssymb,amsfonts}
\usepackage{algorithmic}
\usepackage{graphicx}
\usepackage{textcomp}
\usepackage{bm}
\usepackage{comment}
\usepackage{xcolor}

\usepackage{enumitem}

\newcommand{\norm}[1]{\left\lVert#1\right\rVert}

\newcommand\inner[2]{\left\langle #1, #2 \right\rangle}

\newcommand{\calH}{\mathcal{H}}
\newcommand{\calU}{\mathcal{U}}

\newcommand{\calW}{\mathcal{W}}
\newcommand{\calY}{\mathcal{Y}}
\newcommand{\calX}{\mathcal{X}}
\newcommand{\calZ}{\mathcal{Z}}

\newcommand{\calB}{\mathcal{B}}

\newcommand{\inv}{^{-1}}
\newcommand{\sq}{^{\frac{1}{2}}}

\usepackage{optidef}
\def\BibTeX{{\rm B\kern-.05em{\sc i\kern-.025em b}\kern-.08em
    T\kern-.1667em\lower.7ex\hbox{E}\kern-.125emX}}

\title[Training Lipschitz continuous operators]{Training Lipschitz continuous operators using reproducing kernels}
\usepackage{times}

\author{%
 \Name{Henk J. {van Waarde}} \Email{hvanwaarde@ethz.ch}\\
 \addr Automatic Control Laboratory, ETH Z\"urich, Switzerland
 \AND
 \Name{Rodolphe Sepulchre} \Email{r.sepulchre@eng.cam.ac.uk}\\
 \addr Control Group, University of Cambridge, United Kingdom%
}

\begin{document}

\maketitle

\begin{abstract}%
This paper proposes that Lipschitz continuity is a natural outcome of regularized least squares in kernel-based learning. Lipschitz continuity is an important proxy for robustness of input-output operators. It is also instrumental for guaranteeing closed-loop stability of kernel-based controlllers through small incremental gain arguments. We introduce a new class of \emph{nonexpansive} kernels that are shown to induce Hilbert spaces consisting of only Lipschitz continuous operators. The Lipschitz constant of estimated operators within such Hilbert spaces can be tuned by suitable selection of a regularization parameter. As is typical for kernel-based models, input-output operators are estimated from data by solving tractable systems of linear equations. The approach thus constitutes a promising alternative to Lipschitz-bounded neural networks, that have recently been investigated but are computationally expensive to train. 

\end{abstract}


\section{Introduction}

Although neural networks have impressive performance on tasks such as image recognition, they can be brittle \cite{Szegedy2013}, and are susceptible to adversarial attacks \cite{Eykholt2018}. One way of measuring neural networks' robustness with respect to input perturbations is through their Lipschitz constant (or incremental gain).
To this end, several recent contributions focus on the computation of Lipschitz constants, and the training of neural networks with given Lipschitz constants. Early contributions in this direction are \cite{Szegedy2013,Scaman2018,Tsuzuku2018}. More recently, new methods have emerged that utilize properties of the activation functions: \cite{Combettes2020} treat them as nonexpansive operators while \cite{Fazlyab2019,Pauli2022} exploit incremental quadratic constraints to come up with less conservative bounds. The latter results rely on semidefinite programming and linear matrix inequalities (LMIs). In theory, these LMI-based techniques have polynomial complexity (albeit with high degree). In practice, however, they are computationally demanding and thus only applicable to small-scale networks, c.f. \cite{Revay2020}. 

In this paper we approach the problem of learning Lipschitz continuous input-output operators from a different perspective, namely regularized least squares in reproducing kernel Hilbert spaces (RKHSs). Reproducing kernels play an important role in machine learning algorithms such as support vector machines \cite{Scholkopf2001}. They have also been popularized in systems and control, which led to a renaissance of system identification \cite{Ljung2020}. Starting from the work of \cite{Pillonetto2010}, reproducing kernel Hilbert spaces have proven useful to combine data fitting with prior knowledge of, e.g., stability \cite{Pillonetto2011}, causality \cite{Dinuzzo2015} and positivity \cite{Khosravi2019}. Further contributions focus on reduced order modeling \cite{Bouvrie2017} and error bounds of kernel-based models \cite{Maddalena2021}. 

In contrast to feedforward neural networks that are Lipschitz continuous by the properties of their activation functions\footnote{For example, for nonexpansive activations, a Lipschitz constant is given by the product of the norms of the network's weights \cite{Combettes2020}.}, not all kernel-based models exhibit Lipschitz continuity. This motivates the study of specific classes of kernels. In this paper we introduce so-called \emph{nonexpansive kernels} whose associated Hilbert spaces only contain Lipschitz continuous operators. It is shown that the bilinear, Gaussian a scaled Laplacian kernel are members of this class. We also prove that a Lipschitz constant for any operator in a nonexpansive RKHS is given by the operator's norm. This enables efficient computation of Lipschitz constants for kernel-based models, and paves the way for training operators with given Lipschitz constants via regularized least squares. One of the attractive features of the approach lies in its simplicity. In fact, imposing Lipschitz properties does not lead to computational overhead: identification is performed via standard regularized least squares in suitable (nonexpansive) reproducing kernel Hilbert spaces. 

In addition, we demonstrate how \emph{monotone} operators can be identified using the kernel-based approach, exploiting the Cayley transform. We work with the case that the input and output spaces are general Hilbert spaces, and the results are thus applicable to spaces of square integrable functions and square summable sequences, which are relevant for dynamical systems and control. We illustrate this by identifying a model of the potassium current, one of the components of the Hodgkin-Huxley system describing the behavior of a neuron. 

The outline of the paper is as follows. In Section~\ref{s:background} we provide background material on reproducing kernel Hilbert spaces and regularized least squares. Subsequently, in Section~\ref{s:problem} we state the problem. Section~\ref{s:Lipschitz} contain our results on identifying Lipschitz continuous operators. Then, in Section~\ref{s:monotone} we show how these results can also be applied to train monotone operators. Section~\ref{s:example} treats an illustrative example and Section~\ref{s:conclusions} contains our conclusions. Throughout the paper, we refer to the extended manuscript \cite{vanWaarde2021} for the proofs of the main results. 

\subsection{Notation}

Let $\calX$ and $\calZ$ be real Banach spaces with norms $\norm{\cdot}_\calX$ and $\norm{\cdot}_\calZ$, respectively. We denote the collection of all \emph{bounded linear operators} from $\calX$ to $\calZ$ is denoted by $\calB(\calX,\calZ)$. We denote the operator norm of $A \in \calB(\calX,\calZ)$ by $\norm{A}_{\calB(\calX,\calZ)}$. If $\calX = \calZ$ we simply use the notation $\calB(\calX)$. The identity operator in $\calB(\calX)$ is denoted by $I$. 

Next, let $\calX, \calZ$ be real Hilbert spaces with inner product $\inner{\cdot}{\cdot}_\calX$ and $\inner{\cdot}{\cdot}_\calZ$. We use $A^*$ to denote the \emph{adjoint} of $A \in \calB(\calX,\calZ)$. An operator $A\in\calB(\calX)$ is called \emph{self-adjoint} if $A = A^*$. It is called \emph{positive} if $\inner{x}{Ax}_\calX \geq 0$ for all $x \in \calX$. If $A \in \calB(\calX)$ is a self-adjoint positive operator then there exists a unique self-adjoint positive $A^{\frac{1}{2}} \in \calB(\calX)$, called the \emph{square root}, such that $A = A^{\frac{1}{2}}A^{\frac{1}{2}}$ \cite[p. 265]{Riesz1956}.

\section{Background on reproducing kernel Hilbert spaces}
\label{s:background}

In this section we review some of the theory of reproducing kernel Hilbert spaces of operators \cite{Micchelli2004,Micchelli2005,Carmeli2006,
Caponnetto2008}. 

Throughout, we let $\calU$ and $\calY$ be real Hilbert spaces\footnote{All results remain true if $\calU$ is merely a subset of a Hilbert space.}. In addition, we consider a real Hilbert space $\calH$ of operators from the set of inputs $\calU$ to the output space $\calY$. 

\begin{definition}
\label{d:reproducingkernel}
A mapping $K:\calU\times\calU \to \calB(\calY)$ is called a reproducing kernel for $\calH$ if the following two properties hold:
\begin{itemize}
\item $K(\cdot,u)y : \calU \to \calY$ is a member of $\calH$ for all $u\in\calU$ and $y \in \calY$;
\item The \emph{reproducing property} holds: for all $u\in\calU$, $y \in \calY$ and $H\in\calH$ we have that
\begin{equation}
\label{repprop}
\inner{y}{H(u)}_\calY = \inner{H}{K(\cdot,u)y}_\calH.
\end{equation}
\end{itemize}
We say that $\calH$ is a reproducing kernel Hilbert space if it admits a reproducing kernel. 
\end{definition}

Every reproducing kernel Hilbert space has \emph{exactly one} reproducing kernel \cite[Thm. 1]{Kadri2016}. The class of reproducing kernels is completely characterized by two properties: symmetry and positive semidefiniteness. 

\begin{definition}
A mapping $K: \calU\times \calU \to \calB(\calY)$ is called \begin{itemize}
\item \emph{symmetric} if $K(u,v)^* = K(v,u)$.
\item \emph{positive semidefinite} if for all $n \in \mathbb{N}$, $u_1,u_2,\dots,u_n \in \calU$ and $y_1,y_2,\dots,y_n \in \calY$ we have that
\begin{equation}
\label{defpossemidef}
\sum_{i=1}^n \sum_{j=1}^n \inner{y_i}{K(u_i,u_j)y_j}_\calY \geq 0.
\end{equation} 
\end{itemize}
\end{definition}

The following result, called the Moore-Aronszajn theorem \cite{Aronszajn1950}, shows that $K$ is a reproducing kernel for some Hilbert space if and only if it is symmetric positive semidefinite. For the case that the output is vector-valued, as considered in this paper, we refer to \cite{Micchelli2005}. 

\begin{theorem}
A mapping $K:\calU\times\calU\to\calB(\calY)$ is the reproducing kernel for some reproducing kernel Hilbert space if and only if it is symmetric positive semidefinite. 

Moreover, if $K$ is symmetric positive semidefinite, then there exists a unique reproducing kernel Hilbert space $\calH$ that admits $K$ as a reproducing kernel. 
\end{theorem}

The Moore-Aronszajn theorem is important because it provides a complete classification of reproducing kernel Hilbert spaces: every symmetric positive semidefinite $K$ defines a unique reproducing kernel Hilbert space and vice versa. Some examples of \emph{scalar-valued} reproducing kernels $K :\calU\times\calU \to \mathbb{R}$ are the \emph{polynomial kernel}: 
$$
k(u,v) = (c + \inner{u}{v}_\calU)^d,
$$
with $c \geq 0$ and $d \in \mathbb{N}$, and the \emph{radial basis function kernel}:
\begin{equation}
\label{radialkernel}
k(u,v) = f(\norm{u-v}_\calU^2),
\end{equation} 
where $f: [0,\infty) \to \mathbb{R}$ is completely monotone, i.e., continuous on $[0,\infty)$, infinitely differentiable on $(0,\infty)$ and satisfying
$$
(-1)^n \frac{d^n}{dx^n} f(x) \geq 0
$$
for all $n \in \mathbb{N}$ and $x \in (0,\infty)$ \cite{Schoenberg1938}. The well-known \emph{Gaussian} and \emph{Laplacian} kernels are special cases of radial basis function kernels. More general \emph{vector-valued} $K:\calU\times\calU \to \calB(\calY)$ can be constructed from scalar-valued ones as follows. We have that 
$$
K(u,v) := \sum_{i=1}^r k_i(u,v) R_i
$$ 
is symmetric positive semidefinite if $k_i : \calU\times\calU\to\mathbb{R}$ is symmetric positive semidefinite and $R_i \in \calB(\calY)$ is self-adjoint and positive for all $i = 1,2,\dots,r$ \cite{Kadri2016}. 

Symmetry and positive semidefiniteness can also be expressed in terms of the so-called Gram operator associated with $K$. To do this, we define 
$$
\calY^n := \underbrace{\calY \times \calY \times \cdots \times \calY}_{n \text{ times}}.
$$
Since $\calY$ is Hilbert, $\calY^n$ is a Hilbert space with inner product 
$$
\inner{(y_1,y_2,\dots,y_n)}{(z_1,z_2,\dots,z_n)}_{\calY^n} := \sum_{i = 1}^n \inner{y_i}{z_i}_\calY.
$$
For $u_1,u_2,\dots,u_n \in \calU$ the \emph{Gram operator} $G: \calY^n \to \calY^n$ is defined as 
$$
G(y_1,\dots,y_n) = \left(\sum_{j=1}^n K(u_1,u_j)y_j,\dots,\sum_{j=1}^n K(u_n,u_j)y_j\right).
$$

Then $K$ is symmetric if and only if $G$ is self-adjoint for all $n \in \mathbb{N}$ and $u_1,u_2,\dots,u_n \in \calU$. Moreover, it is positive semidefinite if and only if $G$ is a positive operator for all $n \in \mathbb{N}$ and $u_1,u_2,\dots,u_n \in \calU$, see \cite[Lem. 3]{vanWaarde2021}.

Another way of characterizing reproducing kernels is through the notion of \emph{feature maps}, see also \cite{Micchelli2005}.

\begin{theorem}
\label{t:featuremap}
Let $K : \calU \times \calU \to \calB(\calY)$. Then $K$ is symmetric positive semidefinite if and only if there exists a Hilbert space $\calW$ and a feature map $\phi: \calU \to \calB(\calY,\calW)$ such that for all $u,v \in \calU$: \begin{equation}
\label{featurerep}
K(u,v) = \phi(u)^*\phi(v).
\end{equation}
\end{theorem}

\subsection{Regularized least squares}
\label{s:regLS}

An attractive feature of RKHSs is that several function estimation problems have an elegant and tractable solution if the underlying space has a reproducing kernel. We focus on the regularized least squares problem
\begin{equation}
\label{regLS}
\min_{H \in \calH} \left( \sum_{i=1}^n \norm{y_i - H(u_i)}^2_\calY + \gamma \norm{H}_\calH^2 \right),
\end{equation}
where $y_i \in \calY$, $u_i \in \calU$ for $i = 1,2,\dots,n$ and $\gamma > 0$ is a scalar. It turns out that the solution to \eqref{regLS} is unique for any RKHS, see \cite[Thm. 4.1]{Micchelli2005}. 

\begin{theorem}[RegLS representer theorem]
\label{t:regLS}
Suppose that $\calH$ is a reproducing kernel Hilbert space of operators from $\calU$ to $\calY$ and let $K : \calU \times \calU \to \calB(\calY)$ be its reproducing kernel. There exists a unique solution $\hat{H}$ to \eqref{regLS}, which is given by 
\begin{equation}
\label{estimateregLS}
\hat{H}(\cdot) = \sum_{j=1}^n K(\cdot,u_j) c_j,
\end{equation}
where the coefficients $c_j \in \calY$ ($j = 1,2,\dots,n$) are the unique solution to the system of linear equations
\begin{equation}
\label{lineqregLS}
(G+\gamma I)(c_1,\dots,c_n) = (y_1,\dots,y_n),
\end{equation}
where $G: \calY^n \to \calY^n$ is the Gram operator associated with $K$ and $u_1,u_2,\dots,u_n$.
\end{theorem}

\section{Problem statement}
\label{s:problem}

The goal of this paper is to train robust kernel-based models of the form \eqref{estimateregLS}. One measure of robustness, that has recently been popularized in the literature on neural networks, are bounds on the \emph{Lipschitz constant} \cite{Tsuzuku2018,Fazlyab2019,Pauli2022}. 

\begin{definition}
\label{d:nonexpansive}
Let $\calX$ and $\calZ$ be Banach spaces. An operator $R : D \subseteq \calX \to \calZ$ is 
\begin{itemize}
\item Lipschitz continuous if there exists\footnote{We follow the convention of \cite{Bauschke2011} by considering \emph{a} Lipschitz constant $\ell$ for $R$ and we emphasize that this is not necessarily the smallest $\ell$ satisfying \eqref{Lipschitzcontinuity}.} a nonnegative constant $\ell \in \mathbb{R}$ such that 
\begin{equation}
\label{Lipschitzcontinuity}
\norm{R(x)-R(y)}_\calZ \leq \ell \norm{x-y}_\calX
\end{equation} 
for all $x,y \in D$.
\item Nonexpansive if it is Lipschitz continuous with constant $\ell = 1$, i.e.,
$$
\norm{R(x)-R(y)}_\calZ \leq \norm{x-y}_\calX
$$
for all $x,y \in D$.
\item Contractive if it is Lipschitz continuous with $\ell < 1$.
\end{itemize}
\end{definition}

\noindent
The main problems studied in this paper are to compute Lipschitz constants for kernel-based models, and to train models with \emph{given} Lipschitz constants. We formalize these problems as follows. \\

\noindent 
\textbf{Problem}: Consider $n$ pairs of data points $(u_i,y_i) \in \calU \times \calY$ for $ i = 1,2,\dots,n$.
\begin{enumerate}
\item Let $\calH$ be a reproducing kernel Hilbert space and $\gamma > 0$. Consider the solution $\hat{H}$ to \eqref{regLS}. Find (if it exists) a constant $\ell$ such that \eqref{Lipschitzcontinuity} holds for $\hat{H}$.
\item Given a constant $\ell$, design a suitable reproducing kernel Hilbert space $\calH$ and constant $\gamma > 0$ such that the solution $\hat{H}$ to \eqref{regLS} satisfies \eqref{Lipschitzcontinuity}. 
\end{enumerate}

\section{Identifying Lipschitz continuous operators}
\label{s:Lipschitz}

To start off, we note that for arbitrary RKHSs, the solution to \eqref{regLS} is generally not (globally) Lipschitz continuous. As a simple example, one can consider a polynomial kernel $K(u,v) = (\inner{u}{v}_\calU)^2$. This kernel induces estimates \eqref{estimateregLS} of the form 
$$
\hat{H}(u) = \sum_{j = 1}^n (\inner{u}{u_j}_\calU)^2 c_j,
$$
which do not satisfy \eqref{Lipschitzcontinuity}. This motivates the restriction of the class of symmetric positive semidefinite kernels. In what follows, we recall the definition of \emph{nonexpansive} kernels, introduced in \cite{vanWaarde2021}. This definition is relevant because the reproducing kernel Hilbert space associated to a symmetric positive semidefinite and nonexpansive kernel contains only \emph{Lipschitz continuous} operators. 

\begin{definition}
\label{d:nonexpansivekernel}
A mapping $K: \calU \times \calU \to \calB(\calY)$ is nonexpansive if 
\begin{equation}
\label{nonexpansivekernel}
\norm{K(u,u)-K(u,v)-K(v,u)+K(v,v)}_{\calB(\calY)}^{\frac{1}{2}} \leq \norm{u-v}_\calU.
\end{equation}
holds for all $u,v \in \calU$.
\end{definition}

We emphasize that Definition~\ref{d:nonexpansivekernel} introduces nonexpansiveness for \emph{two-variable} mappings, and is thus different from the classical Definition~\ref{d:nonexpansive}. Nonetheless, the term ``nonexpansive" is natural also in the context of kernels because it is intimately related to Definition~\ref{d:nonexpansive} through \emph{feature maps}. Indeed, let $K:\calU\times\calU\to\calB(\calY)$ be symmetric positive semidefinite. Then, by Theorem~\ref{t:featuremap}, $K(u,v) = \phi(u)^* \phi(v)$ for some $\phi: \calU \to \calB(\calY,\calW)$, where $\calW$ is a Hilbert space. Note that
\begin{equation*}
\begin{aligned}
K(u,u)-K(u,v)-K(v,u)+K(v,v) = (\phi(u) - \phi(v))^*(\phi(u) - \phi(v))&.
\end{aligned}
\end{equation*}
By \cite[Fact 2.18(ii)]{Bauschke2011}, this means that
\begin{align*}
\norm{K(u,u)-K(u,v)-K(v,u)+K(v,v)}_{\calB(\calY)}^{\frac{1}{2}} = \norm{\phi(u) - \phi(v)}_{\calB(\calY,\calW)}&.
\end{align*}
Therefore, a kernel $K$ is nonexpansive if and only if all its associated feature maps are nonexpansive in the sense of Definition~\ref{d:nonexpansive}. 

\begin{theorem}
\label{t:nonexpansive}
Let $\bar{u} := (u_1,u_2,\dots,u_n) \in \calU^n$ and $\bar{y} := (y_1,y_2,\dots,y_n) \in \calY^n$ be data. 
Consider a symmetric positive semidefinite kernel $K: \calU\times \calU \to \calB(\calY)$, and let $\calH$ be its associated reproducing kernel Hilbert space. Assume that $K$ is nonexpansive. Then the following statements hold:
\begin{enumerate}[label=(\alph*)]
\item Every $H \in \calH$ is Lipschitz continuous with constant $\norm{H}_\calH$. \label{t:nonexpansive1}
\item The solution $\hat{H} \in \calH$ to the regularized least squares problem \eqref{regLS} has norm
$$
\norm{\hat{H}}_\calH = \norm{G\sq(G+\gamma I)\inv \bar{y}}_{\calY^n},
$$
where $G :\calY^n \to \calY^n$ is the Gram operator associated with $K$ and $\bar{u}$.
\label{t:nonexpansivenormHhat}
\item Thus, $\hat{H}$ is Lipschitz continuous with constant $\ell$ if $\gamma > 0$ satisfies
\begin{equation}
\label{conditiongamma}
\norm{G\sq(G+\gamma I)\inv \bar{y}}_{\calY^n} \leq \ell. 
\end{equation}
\label{t:nonexpansive2}
\end{enumerate}  
\end{theorem}

A proof of Theorem~\ref{t:nonexpansive} is provided in \cite{vanWaarde2021}. We can now draw a few conclusions: 
\begin{enumerate}
\item In reproducing kernel Hilbert spaces associated with nonexpansive kernels, every operator is Lipschitz continuous, and a Lipschitz constant is given by the RKHS norm of the operator. 
\item To train robust kernel-based models satisfying \eqref{Lipschitzcontinuity}, we can apply regularized least squares with i) a suitable (nonexpansive) kernel, and ii) a suitable parameter $\gamma$ satisfying \eqref{conditiongamma}. 
\item An attractive feature of this approach is its simplicity. In fact, by Theorems~\ref{t:regLS} and \ref{t:nonexpansive}, Lipschitz continuous operators can be identified by solving tractable systems of linear equations. This is a potential advantage over training robust neural networks \cite{Fazlyab2019,Pauli2022} that relies on semidefinite programming. 
\end{enumerate}

In the following proposition from \cite{vanWaarde2021}, we highlight some examples of nonexpansive kernels.

\begin{proposition}
\label{p:nonexpansivekernels}
The following (scalar-valued) symmetric positive semidefinite kernels $k: \calU \times \calU \to \mathbb{R}$ are nonexpansive:
\begin{itemize}
\item The bilinear kernel $\inner{u}{v}_\calU$.
\item The Gaussian kernel $e^{\frac{-\norm{u-v}^2_\calU}{\sigma^2}}$ whenever $\sigma \geq \sqrt{2}$. 
\item The scaled Laplacian kernel 
\begin{equation*}
k(u,v) = (1+\norm{u-v}_\calU) e^{-\norm{u-v}_\calU}.
\end{equation*}
\item The kernel
\begin{equation*}
k(u,v) = (c + \norm{u-v}_\calU^2)^{-d}
\end{equation*} 
where $c \geq 0$ and $d > 0$ are reals satisfying $2d \leq c^{d+1}$. 
\end{itemize} 

In addition, $K: \calU \times \calU \to \calB(\calY)$ is both symmetric positive semidefinite and nonexpansive if
\begin{itemize}
\item $K(u,v) = k(u,v) R$ where $k : \calU \times \calU \to \mathbb{R}$ is a symmetric positive semidefinite and nonexpansive kernel, and $R \in \calB(\calY)$ is self-adjoint and positive with $\norm{R}_{\calB(\calY)} \leq 1$. 
\item $K(u,v) = \sum_{i = 1}^r \alpha_i K_i(u,v)$, where $K_i : \calU \times \calU \to \calB(\calY)$ is nonexpansive symmetric positive semidefinite, $\alpha_i \geq 0$ for all $i = 1,2\dots,r$ and $\sum_{i = 1}^r \alpha_i \leq 1$. 
\item $K(u,v) = R L(u,v) R^*$, where $L:\calU \times \calU \to \mathbb{R}$ is a symmetric positive semidefinite and nonexpansive kernel, and $R \in \calB(\calY)$ satisfies $\norm{R}_{\calB(\calY)} \leq 1$. 
\end{itemize}
\end{proposition}

\section{Monotone operators}
\label{s:monotone}

A concept that is closely related to nonexpansiveness is \emph{monotonicity}.

\begin{definition}
\label{d:monotone}
Let $\calX$ be a Hilbert space and $R : D \subseteq \calX \to \calX$. We say that $R$ is monotone if
$$
\inner{x-y}{R(x)-R(y)}_\calX \geq 0
$$
for all $x,y \in D$.
\end{definition}

Monotone operators play a fundamental role in convex analysis and optimization \cite{Ryu2022}. In the special case that $\calX$ is the space of square integrable functions, monotonicity is also closely related to \emph{incremental passivity} (and the notions coincide for causal operators \cite{Desoer1975}). Monotonicity can thus be interpreted as the physical property that any input/output trajectory $(x,R(x))$ can only dissipate energy \emph{with respect to} any other given trajectory $(y,R(y))$.

The following proposition asserts a well-known relation between monotone and nonexpansive operators \cite{Bauschke2011}. 

\begin{proposition}
\label{p:scattering}
Let $\calX$ be a Hilbert space and $R: \calX \to \calX$. Suppose that $I+R$ is invertible. Define the operator $S: \calX \to \calX$ as
\begin{equation}
\label{operatorS}
S := (I-R)(I+R)\inv.
\end{equation}
Let $(u,y),(v,z) \in \calX \times \calX$ be related by
$v = u+y$ and $z = u-y$. 
Then we have that
\begin{enumerate}[label=(\roman*)]
\item \label{i:MN} $I+S$ is invertible and
\begin{equation}
\label{operatorR}
R = (I-S)(I+S)\inv.
\end{equation} 
\item \label{i:yTu} $y = R(u)$ if and only if $z = S(v)$.
\item \label{i:TsPhi} $R$ is monotone if and only if $S$ is nonexpansive.
\end{enumerate}
\end{proposition}

\noindent
The operation on $S$ in \eqref{operatorR} is often referred to as the \emph{Cayley transform}. The consequence of Proposition~\ref{p:scattering} is that we can also identify monotone operators using reproducing kernel Hilbert spaces. To this end, let us assume that $\calU = \calY$. Given data samples $u_1,u_2,\dots,u_n \in \calU$ and $y_1,y_2,\dots,y_n \in \calY$, we can compute the ``transformed data" $v_i := u_i + y_i$ and $z_i = u_i - y_i$, and apply Theorem~\ref{t:nonexpansive} (with $\ell = 1$) to identify a nonexpansive operator $S$ from the input/output data $(v_i,z_i)$ for $i = 1,2,\dots,n$. The operator $R$, defined in terms of $S$ in \eqref{operatorR}, is then monotone by Proposition~\ref{p:scattering}. 

In some cases, it is beneficial to choose a slightly larger regularization parameter in Theorem~\ref{t:nonexpansive}, so that $\ell$ is strictly less than $1$ and $S$ is a contraction. Indeed, in this case, the input/output behavior of $R$ can be simulated efficiently using fixed point algorithms. 

\begin{proposition}
\label{p:fixedpoint}
Suppose that $S:\calU \to \calY$ is a contraction with Lipschitz constant $\ell \in [0,1)$. Then $I+S$ is invertible so $R$ in \eqref{operatorR} is well-defined. Moreover, for any $u^* \in \calU$, the output $y^* = R(u^*)$ can be computed via the Picard iteration $y^* = \lim_{k \to \infty} y^k$, where $y^0 \in \calY$ is arbitrary and
\begin{equation*}
y^{k+1} = u^* - S(u^* + y^k) \text{ for } k \geq 0.
\end{equation*} 
In addition, for any $k \geq 0$ we have that 
\begin{equation*}
\norm{y^k - y^*}_\calY \leq \ell^k \norm{y^0 - y^*}_\calY.
\end{equation*} 
\end{proposition}

The proposition follows from the Banach fixed point theorem and the fact that the mapping $y \mapsto u^*-S(u^*+y)$ is a contraction by the hypothesis on $S$. 

\section{Illustrative example}
\label{s:example}

We consider a model of the potassium ion channel, which is a component of the Hodgkin-Huxley dynamical system \cite{Hodgkin1952} describing the electrical characteristics of an excitable cell. This \emph{conductance-based} model relates the current through the potassium channel and the voltage across the cell's membrane through a nonlinear conductance. It is described by
\begin{equation}
\begin{aligned}
\label{potassiumconductance}
\dot{x} &= \alpha(u)(1-x) - \beta(u)x, \:\: x(0) = 0 \\
y &= g x^4 (u-\bar{u}),
\end{aligned}
\end{equation}
where $x: \mathbb{R} \to \mathbb{R}$ denotes the gating (state) variable, $u: \mathbb{R} \to \mathbb{R}$ denotes the (input) voltage across the cell's membrane, and $y : \mathbb{R} \to \mathbb{R}$ is the (output) potassium current per unit area. The functions $\alpha$ and $\beta$ depend on the voltage but not explicitly on time. Both $\bar{u} \in \mathbb{R}$ and $g \in \mathbb{R}$ are constants, called the \emph{reversal potential} and \emph{maximal conductance}, respectively. 

Owing to the original proposal by Hodgkin and Huxley, the constants $g$ and $\bar{u}$ are chosen as $g = 36$ and $\bar{u} = 12$, while the functions $\alpha$ and $\beta$ are given by
\begin{align*}
\alpha(u) &= 0.01 \frac{u+10}{e^{\frac{u+10}{10}}-1} \\
\beta(u) &= 0.125 e^{\frac{u}{80}}.
\end{align*}

\begin{figure}[h!]
\centering
\includegraphics[width=.9\textwidth]{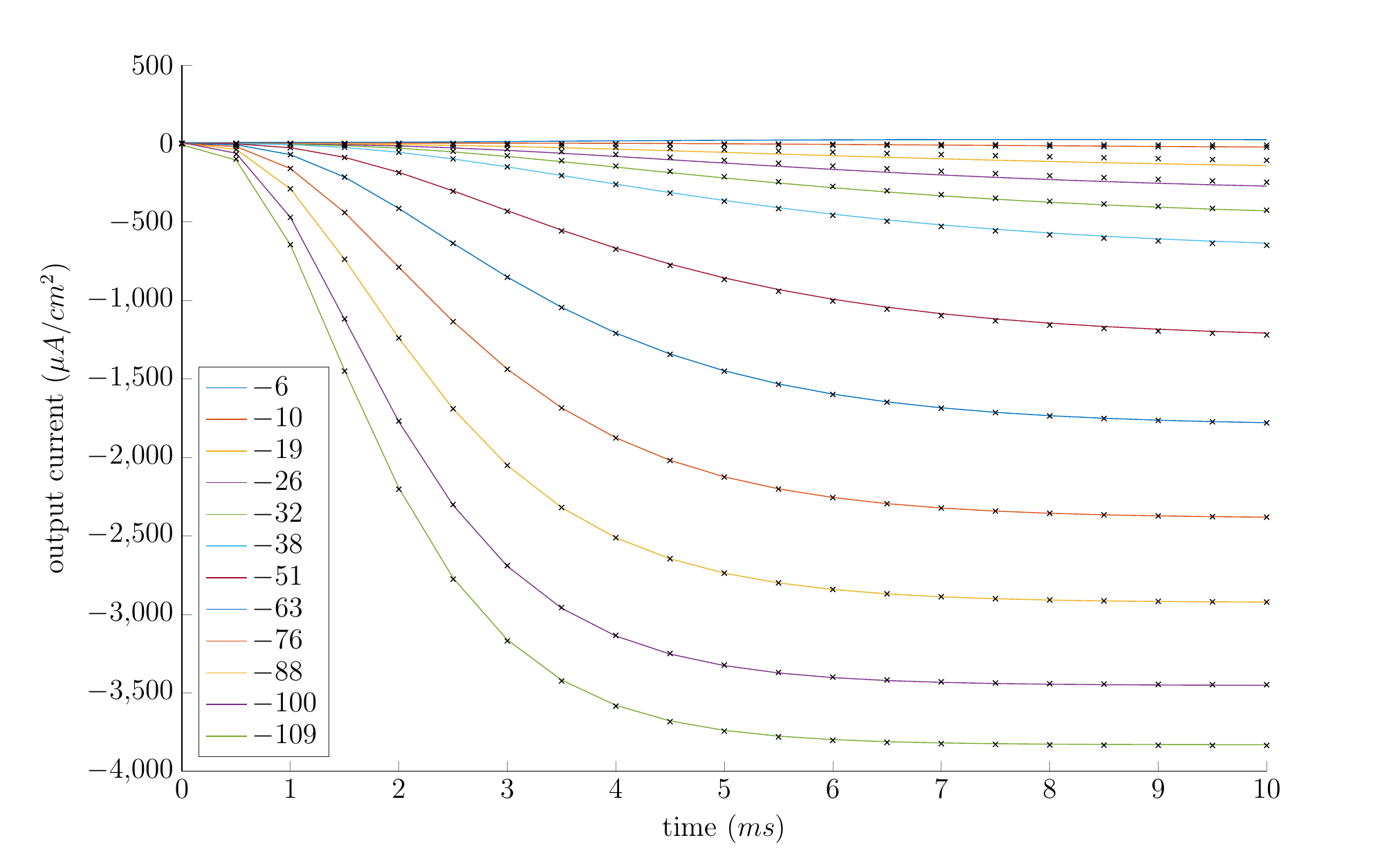}
\caption{Outputs of the identified operator for constant input voltages from $-6 mV$ to $-109 mV$ (in colors). The black crosses correspond to the data samples at times $t_0,t_1,\dots,t_{20}$.}
\label{fig:identifiedoperator}
\end{figure}

Hodgkin and Huxley determined the model \eqref{potassiumconductance} and its parameters on the basis of a series of step response experiments where the voltage $u$ was kept at different constant values, ranging from $-6mV$ to $-109mV$, and the corresponding current $y$ was measured, see \cite[Part II]{Hodgkin1952}. Using the model \eqref{potassiumconductance}, we have reconstructed the input-output data of \cite{Hodgkin1952} over the time instants $t_0,t_1,\dots,t_{20}$, where $t_j = 0.5j$, and the results are displayed as black crosses in Figure~\ref{fig:identifiedoperator}. 

All these experimental data suggest that the input-output behavior of the potassium current defines a \emph{monotone} operator on $\ell_2(\{0,1,\dots,20\})$. Nonetheless, it can be proven that the model \eqref{potassiumconductance} is \emph{not} monotone \cite{vanWaarde2021}. Next, we will identify a monotone operator of the potassium current using the tools of this paper. We use a kernel of the form
$$
K(u,v) = (1+\norm{u-v}_{\ell_2}) e^{-\norm{u-v}_{\ell_2}} \cdot I,
$$
which is nonexpansive by Proposition~\ref{p:nonexpansivekernels}. We choose a value of $\gamma = 4.441\cdot 10^{-4}$ which implies that the left hand side of \eqref{conditiongamma} is $0.9903 < 1$. This results in an identified nonexpansive operator of the form \eqref{estimateregLS}. Lastly, we exploit the fixed point algorithm in Proposition~\ref{p:fixedpoint} to simulate the Cayley transform of this identified system. The simulation results are reported in Figure~\ref{fig:identifiedoperator} for different constant input values (in various colors). Note that the curves are obtained by interpolating between the output values at times $t_j$ ($j = 0,1,\dots,20$).  

We observe that the identified operator explains the data well, with a small misfit for larger input values and near-perfect reconstruction for smaller ones. Importantly, by Proposition~\ref{p:scattering} the identified operator is monotone.

\newpage 

\section{Conclusions}
\label{s:conclusions}

In this paper, we have introduced a method to incorporate bounds on the Lipschitz constant in kernel-based regularized least squares problems. As our main result, we have introduced a new class of \emph{nonexpansive} kernels that induce Hilbert spaces consisting of Lipschitz continuous operators. Using regularized least squares, Lipschitz continuous operators can then be identified from data and their Lipschitz constants can be tuned by appropriate choice of the regularization parameter. We have also demonstrated how the approach enables the identification of monotone operators via the Cayley transform. The limitation of state-space modeling to achieve this objective was illustrated with a simple model of a nonlinear circuit, in particular the celebrated model of the potassium current of \cite{Hodgkin1952}.

In the context of machine learning, the results of this paper are well-aligned with the objectives of \cite{Tsuzuku2018,Fazlyab2019,Pauli2022}. These papers focus on training input-output operators satisfying Lipschitz properties using feedforward neural network architectures. In this paper, we have argued that kernel-based models are easier to train, but more studies are needed to compare the performance and to understand the relative merits of both approaches. 

In the context of systems and control, Lipschitz continuity and monotonicity are closely related to the notions of \emph{finite incremental gain} and \emph{incremental passivity} \cite{Desoer1975,vanderSchaft2017}. These properties are the cornerstone of the analysis of feedback systems. This relation is explored in more detail in \cite{vanWaarde2021}. 

\bibliography{references}

\end{document}